\documentclass[12pt]{amsart}
\usepackage{txfonts}
\usepackage{amscd}
\usepackage{cite}
\usepackage{epsfig}
\usepackage{verbatim}
\usepackage{mathdots}
\usepackage{amssymb}
\usepackage{amsfonts}
\usepackage{amsbsy}
 \usepackage{graphicx}
 \usepackage[usenames]{color} 
 \usepackage{cancel}

\newtheorem{theo}{Theorem}[section]
\newtheorem{lem}[theo]{Lemma}
\newtheorem{prop}[theo]{Proposition}
\newtheorem{coro}[theo]{Corollary}

\theoremstyle{definition}
\newtheorem{defi}[theo]{Definition}

\numberwithin{equation}{section}

\newcommand{\R}{{\mathbb R}}
\newcommand{\Z}{{\mathbb Z}}

\newcommand{\Q}{{\mathbb Q}}
\newcommand{\N}{{\mathbb N}}

\newcommand{\C}{{\mathcal C}}
\newcommand{\D}{{\mathcal D}}
\newcommand{\E}{{\mathcal E}}

\newcommand{\eproof}{\hfill$\square$}

\begin{document}

\title{On spectral $N$-Bernoulli measures }

\author{Xin-Rong Dai}
\address{School of Mathematics and Computational Science, Sun Yat-Sen University, Guangzhou, 510275,  P. R. China}
\email{daixr@mail.sysu.edu.cn}
\author[X.-G. He]{Xing-Gang He}
\address{School of Mathematics and Statistics   \\ Central China Normal University\\ Wuhan, 430079
\\P. R. China.}
\email{xingganghe@163.com}
\author[K. -S. Lau]{Ka-Sing Lau}
\address{Department  of Mathematics   \\ The Chinese University of Hong Kong\\Hong Kong}
\email{kslau@math.cuhk.edu.hk}
\thanks{The research is partially supported by the RGC grant of Hong Kong; The second author and the third author are also supported by
the NNSF of
China 11271148  and 11171100 respectively. 
}

\maketitle

\begin{abstract}\ For $0<\rho<1$ and  $N>1$ an integer, let $\mu$ be the
self-similar measure defined by $\mu(\cdot)=\sum_{i=0}^{N-1}\frac
1N\mu(\rho^{-1}(\cdot)-i)$. We prove that $L^2(\mu)$ has an exponential
orthonormal basis if and only if $\rho=\frac 1q$ for some $q>0$ and $N$ divides
$q$. The special case is the Cantor measure with $\rho =\frac 1{2k}$
and $N=2$ \cite {JP}, which was proved recently to be the only
spectral measure among the Bernoulli convolutions with $0<\rho<1$
\cite {D}.
\end{abstract}

\tableofcontents

\bigskip

\section{\bf Introduction}
\setcounter{equation}{0}

Let $\mu$ be a  probability measure on ${\Bbb R}^s$ with compact
support.
  For a countable subset $\Lambda \subset {\Bbb R}^s$, we let   $e_\Lambda=\{e_\lambda=e^{-2\pi i \langle\lambda, x\rangle}:
\lambda\in\Lambda\}.$  We call $\mu$ a spectral measure, and  $\Lambda$  a {\it spectrum}  of $\mu$
if $e_\Lambda$ is an orthogonal basis for $L^2(\mu)$. The existence and nonexistence of a spectrum for $\mu$ is a  basic problem in harmonic
analysis, it was initiated by Fuglede  in his seminal paper {\cite F},  and has been studied extensively since then \cite{D, DHL, DHS,
DHSW, DJ, DL, HLL, JP, LW, Lai, Lan, Li1, Li2, S, T}. Recently
He, Lai and Lau \cite{HLL} proved that a spectral measure $\mu$ must
be of pure type, i.e., $\mu$ is absolutely continuous or singular
continuous with respect to the Lebesgue measure or counting measure
supported on a finite set (actually this holds more generally for
{\it frames}). When $\mu$ is the Lebesgue measure restricted on a
set $K$ in $\R^s$, it is well-known that the spectral property is
closely connected with the tiling property of $K$, and is known as
the Fuglede problem \cite {{F}, {LW}, {KM}, {T}}. For continuous
singular measures, the first spectral measure was given by Jorgensen
and Pedersen \cite{JP}: the Cantor measure $\mu_\rho$ with
contraction ratio $\rho =1/2k$. There are considerable studies for
such measures \cite {{DHL}, {DHS}, {HuL}, {LW}, {Li1}, {Li2},
{S}}, and a celebrated open problem was to characterize the
spectral measures $\mu_\rho, 0<\rho <1$ among the Bernoulli
convolutions
$$
\mu_\rho(\cdot)=1/2\mu_\rho(\rho^{-1}\cdot) + 1/2\mu_\rho(\rho^{-1}\cdot -1).
$$
In \cite {HuL}, Hu and Lau showed that $\mu_\rho$ admits an infinite
orthonormal set if and only if $\rho$ is the $n$-th root of $p/q$
where $p$ is odd and $q$ is even. The characterization problem was
finally completed recently by Dai \cite {D} that the above Cantor
measures $\mu_{1/(2k)}$ is the only class of spectral measures among
the $\mu_\rho$.

\bigskip
In this paper we study the spectrality of the self-similar measures.
Let $0< \rho <1$,  $\D =\{0, d_1, \cdots , d_{N-1}\}$ a finite set
in $\R$, and $\{w_j\}_{j=0}^{N-1}$ a set of probability weights. We
call $\mu$ a {\it self-similar measure}  generated by $(\rho, \D)$
and $\{w_j\}_{j=0}^{N-1}$ if $\mu$ is the unique probability measure
satisfying
\begin{equation} \label {eq1.1}
\mu(\cdot )=\frac 1{N}\sum_{j=0}^{N-1} w_j \ \mu(\rho^{-1} (\cdot )
-d_j).
\end{equation}
We will use $\mu_{\rho, N}$ to denote the special case where  $\D =\{0, \cdots , N-1\}$ with uniform weight, i.e.,
\begin {equation} \label {eq1.2}
\mu_{\rho, N}(\cdot )=\frac 1N \sum_{j=0}^{N-1}\
\mu_{\rho, N}(\rho^{-1}(\cdot ) - j).
\end{equation}
The spectral property of such measure was first studied  by Dai, He
and Lai \cite {DHL} as a generalization of the Bernoulli convolution in \cite {D} (${\mathcal D}=\{0,1\}$).
Our main result in this paper is to extend the characterization of spectral Bernoulli convolution to the class of  $\mu_{\rho, N}$ in (\ref{eq1.2}).
 Our motivation to extend the Bernoulli convolutions to this class of measures is due to a conjecture of $\L$aba and Wang, and also to answer a
 question on the convolution of spectral measures (see the remark in \S 6).
We prove

\medskip

\begin{theo} \label{th1.1} Let  $0<\rho<1$.  Then $\mu_{\rho, N}$ is a spectral measure
if and only if $\rho=\frac 1q$ for some integer $q>1$ and $N\mid q$.
\end{theo}

\medskip
The sufficiency of the theorem follows from the same pattern as the Cantor measure in \cite{JP} by producing a Hardamard matrix, then construct
 the  canonical spectrum (see Section 2). Ont the other hand, the proof of the necessity needs more work. We observe that for $\mu_{\rho, N}$ to be a spectral measure,
 $\rho$ must be an algebraic number. we prove by elimination that each of the following cases can NOT admit a spectrum: (for $\frac pq$,
 we always assume they have no common factor)

\vspace{0.1cm}

\ \ \ (i) \ $\rho=(\frac pq)^{1/r}$ for some $r>1$ (it is an irrational),
\ (Proposition \ref{pro2.3});

\vspace{0.1cm}

\ \ (ii) \ $\rho \not = (\frac pq)^{1/r}$ for any $r>1$ and is an irrational (Proposition \ref{pro2.6});

\vspace{0.1cm}

\  (iii) \ $\rho=\frac pq $ and $1\le \hbox {gcd}(N, q) <N$ \    (Proposition \ref{prop3.1});

\vspace{0.1cm} \  \ (iv) \ $\rho=\frac pq $, $p>1$ and $N\mid q$ \
(Proposition \ref{prop3.2}).

\medskip

Let $\widehat \mu_{\rho, N}$ be the Fourier transform of $\mu_{\rho, N}$, and ${\mathcal Z}(\widehat \mu_{\rho, N})$ the zeros
of $\widehat \mu_{\rho, N}$.  The proof is based on the criteria in Theorem \ref {th2.1} and Lemma \ref{th2.2}, and the technique
 is to make use of some explicit expressions  of ${\mathcal Z}(\widehat \mu_{\rho, N})$, and that
 $\Lambda - \Lambda \subset {\mathcal Z}(\widehat \mu_{\rho, N})$ for any exponential orthogonal set $\Lambda$.

\medskip

The most subtle part of the proof is (iv). As is
known, there is certain canonical  $q$-adic expansion
of $\lambda$ in a spectrum $\Lambda$ (see (\ref{eq3.1})), and
there are also  others.  In \cite {DHS}, Dutkay {\it et al}
treated the $4$-adic expansions as in a symbolic space $\Omega_2^*$,
and consider certain maps on $\Omega_2^*$ to $\Z^+$ to preserve the
maximal orthogonality property.
 This idea was refined and investigated by Dai, He and Lai \cite {DHL} by
 replacing the $q$-adic expansion on $\Z$  with digits in $ \C= \{-1,0,\cdots, q-2\}$. Let $\iota : \Omega_N^* \to \C$ be a {\it selection
map} as defined in Definition \ref{def4.1} (it was called a {\it
maximal map} in \cite{DHL}), and let $\iota^* ({\bf i})
=\sum_{n=1}^\infty \iota({\bf i}0^\infty|_n) q^{n-1}$. The importance
of the selection map is in the following theorem (Theorem \ref{thm3.5}), which also has independent interest.

\medskip

\begin{theo} \label{th1.2}  Suppose $\rho =p/q$ and $N\mid q$. Then $\Lambda \subset {\mathcal Z}(\widehat{\mu}_{\rho, N})$
defines a maximal exponential orthogonal subset in $L^2({\mu}_{\rho, N})$  if and only if there exist $m_0\geq 1$ and a
selection map $\iota$ such that
$\Lambda=\rho^{-m_0}N^{-1}(\iota^*(\Omega^{\iota}_N))$.
\end{theo}

 \medskip

We organize the paper as follows. In Section 2, we set up the notations, the basic criteria of spectrum, and the element properties of the zero set
 ${\mathcal Z}(\mu_{\rho, N})$. We settle cases (i), (ii) in Section 3. For the case $\rho = p/q$, in Section 4 we give a detail study
 of the  maximality of $\Lambda$ such the $\Lambda-\Lambda \subset {\mathcal Z}(\mu_{\rho, N})$,  which is used in Section 5 to consider
  cases (iii) and (iv) . In Section 6, we give some remarks of the spectral measures and the  remaining questions.

\bigskip
\bigskip
\section {\bf Preliminaries}

We assume that  $\mu$ is a probability measure with compact support.
The Fourier transformation of $\mu$ is define as usual,
$$
\widehat{\mu}(\xi)=\int e^{-2\pi i \xi x}d\mu(x).
$$
Let ${\mathcal Z}(\widehat{\mu}):=\big \{\xi:
\widehat{\mu}(\xi)=0\}$ be the set of zeros of $\widehat{\mu}$.  We
denote the complex exponential function $e^{-2 \pi i \lambda (\cdot
)}$ by $e_{\lambda}$. Note that $\{e_\lambda: \lambda\in\Lambda\}$
is an orthogonal set in $L^2(\mu)$ if and only if $\widehat{\mu}
(\lambda_i-\lambda_j) =0$ for any $\lambda_i \ne \lambda_j \in
\Lambda$; $\Lambda$ is called a {\it spectrum} of $ \mu$ if
$\{e_\lambda\}_{\lambda \in \Lambda}$ is an orthonormal basis for
$L^2(\mu)$. For $\xi \in {\Bbb R}$, we let
 $$
 Q(\xi)=\sum_{\lambda\in\Lambda}|\widehat{\mu}(\xi+\lambda)|^2.
 $$
  The following  theorem is a basic criterion for the spectrality of $\mu$ \cite{JP}.

\medskip
\begin{theo} \label{th2.1} Let $\mu$ be a probability measure with compact support, and let $\Lambda\subset \R$ be a countable subset. Then

\vspace {0.2cm}
\ (i) \ $\{e_\lambda\}_ {\lambda \in \Lambda}$ is an orthonormal set  of $L^2(\mu)$ if
and only if $Q(\xi)\le 1$ for $\xi\in\R$; and

 \vspace {0.1cm}
(ii)\ \ it is an orthonormal basis if and
only if $Q(\xi)\equiv 1$ for $\xi\in\R$.
\end{theo}

\medskip

Throughout the paper, we use the notation $\Lambda$ to denote a
subset such that $0 \in \Lambda$ and $\Lambda \setminus \{0\}
\subset {\mathcal Z}(\widehat \mu)$. We say that $\Lambda$ is a {\it
bi-zero set} of $\mu$ if $(\Lambda - \Lambda)\setminus \{0\} \subset
{\mathcal Z}(\widehat{\mu})$, and call it a {\it maximal} bi-zero
set if it is  maximal in ${\mathcal Z}(\widehat{\mu})$ to have the
set difference property.
 Clearly that $\Lambda$ is a bi-zero set is equivalent to $\{e_\lambda: \lambda \in \Lambda  \}$ is an orthogonal subset of  $L^2(\mu)$. An exponential
 orthonormal basis corresponds to a maximal
bi-zero set, but the converse is not true. In fact  we will give a
characterize of the maximal bi-zero sets of $\mu_{\rho, N}$
for the case $\rho=\frac pq$ and $N\mid q$  in Section
4, and establish the spectrality through Theorem \ref{th2.1}(ii) in
Section 5.

 \bigskip

As a simple consequence of Theorem \ref{th2.1}, we have the following useful lemma.

\medskip

\begin{lem}\label{th2.2} Let $\mu=\mu_0\ast \mu_1$ be the convolution of two probability measures $\mu_i$, $i=0, 1$,  and they are not
Dirac measures. Suppose that $\Lambda$ is a bi-zero set of $\mu_0$,
then $\Lambda$ is also a bi-zero of $\mu$, but cannot be a spectrum
of $\mu$.
\end{lem}

\medskip

\noindent {\bf Proof.} Note that $\mu_i$ is not an Dirac measure is
equivalent to  $|\widehat{\mu}_i(\xi)|\nequiv 1$. Since
$\widehat{\mu}_0(0)=1$, there exists  $\xi_0$ such that
$|\widehat{\mu}_0(\xi_0)|\ne 0$ and $|\widehat{\mu}_1(\xi_0)|<1$.
Hence by Theorem \ref{th2.1}(i),
$$
Q(\xi_0)=\sum_{\lambda\in\Lambda}|\widehat{\mu}(\xi_0+\lambda)|^2=\sum_{\lambda\in\Lambda}
|\widehat{\mu}_0(\xi_0+\lambda)|^2|\widehat{\mu}_1(\xi_0+\lambda)|^2<\sum_{\lambda\in\Lambda}
|\widehat{\mu}_0(\xi_0+\lambda)|^2\le 1.
$$
The result follows by Theorem \ref{th2.1}(i) and (ii).\eproof

\bigskip

Now we consider the self-similar measure $\mu_{\rho, N}$ in Theorem \ref{th1.1}. It was proved in \cite{DHL} that if $\rho = 1/q$ and $N\mid q$,
 then $\mu_{\rho, N}$ is a spectral measure.  The proof is quite simple. In fact as $N\mid q$, we write $q=Nr$.  If $r=1$, then $\mu$ is just the
 Lebesgue measure on the unit interval, and the result is trivial. If $r >1$,
 observe that  for $\D=\{0, \cdots , N-1\}$ and $\Gamma =r\{0,
\cdots, N-1\}$,  the matrix
$$
H:=[e^{2\pi i \frac {i k}{q}}]_{i \in \D, k \in \Gamma} =[e^{2\pi i
\frac {i j}{N}}]_{0\leq i,j \leq N}
$$
is a Hadamard matrix (i.e., $HH^* =NI$). This shows that $({q}^{-1} \D, \Gamma)$ is a {\it compatible pair},
 hence $\mu_{1/q, N}$ is a spectral measure \cite {LW},
and the canonical spectrum is given by
\begin {equation}\label {eq3.1}
\Lambda =\big  \{\sum_{j=0}^k a_j q^j: \ a_j \in \Gamma, \ k\geq
0\big\}
\end{equation}
(note that the spectrum is not unique). Our main task is to prove the converse. The strategy is to eliminate all the possible cases
so that the only admissible case is $\rho =1/q$ with $N\mid q$.

\bigskip

Recall that the Fourier transform of $\mu_{\rho, N}$ has the following expression
$$
\widehat{\mu}_{\rho, N}(\xi)=M_N(\rho\xi)\ \widehat{\mu}_{\rho,
N}(\rho\xi)=\prod_{k=1}^\infty M_ N(\rho^k\xi)
$$
where $ M_N(\xi)=\frac 1N\sum_{j=0}^{N-1} e^{-2\pi i j\xi} $ is the
{\it mask polynomial} of $\D$. It is clear that  $
|M_N(\xi)|=\left|\frac{\sin N\pi \xi}{N\sin \pi\xi}\right|, $ and
the zeros  of $M_N(\xi)$ is $a/N, \  a \in {\Bbb Z} \setminus \{0\},
\ N\nmid a$.\ Let
\begin{equation} \label {eq2.0}
{\mathcal Z}(M_N)= \Big \{\frac aN :  a \in {\Bbb Z} \setminus
\{0\}, N \not \mid a \Big \} = \Big \{\frac aN:  a\in\Z \setminus
N\Z \Big \}.
\end {equation}
It  follows from the infinite product expression of $\widehat \mu_{\rho, N}$ that
\begin{equation} \label {eq2.1}
{\mathcal Z}(\widehat{\mu}_{ \rho, N})=\big \{\rho^{-k}\frac aN:\
k\ge 1, \  a\in\Z \setminus N\Z \big \}.
\end{equation}

\bigskip
 For distinct $\lambda_1,\lambda_2 \in \Lambda \setminus \{0\}$, (\ref{eq2.1}) and the bi-zero property of $\Lambda$ imply that
$$
\rho^{-k_1} \frac {a_1}N - \rho^{-k_2} \frac {a_2}N = \rho^ {-k} \frac aN.
$$
Hence $\rho$ is an algebraic number.  Recall that an {\it algebraic
number} is a root of an integer equation of the form $c_0 x^n+ c_1
x^{n-1} +\cdots + c_n \in {\Bbb Z}[x]$, and it is called an {\it
algebraic integer} if $c_0=1$.

\bigskip
\bigskip

\section {\bf Spectrality for irrational contraction }

For any integer $r\geq 1$, let
$$
{\Bbb Q}^{1/r} = \big \{ \rho = u^{1/r} : \  0< u <1 \ \ \hbox {is a rational}\big \}.
$$
We make the convention that the above $r$ is the smallest integer for $\rho = u^{1/r}$ (for example:  $\rho = (\frac 49)^{1/4} = (\frac 23)^{1/2}$, we will take $r =2$). Hence for $\rho\in {\Bbb Q}^{1/r}, r>1$, then $\rho$ is an irrational.

\medskip

\begin{prop}\label{pro2.3} Let $\rho\in {\Bbb Q}^{1/r}, r>1$, then $\mu_{\rho, N}$ is not a spectral measure.
\end{prop}

\medskip

\proof Let $\rho = u^{1/r}$ where $0<u<1$ is a rational.  We write
$$
\widehat{\mu}_{\rho, N}(\xi)=\prod_{k=1}^\infty
M_N(\rho^k\xi)=\prod_{k=0}^\infty \prod_{i=1}^{r} M_N\big
(u^k\rho^i\xi\big ).
$$
Define the probability measures $\mu_i(\cdot ) = \mu_{u,
N}(u\rho^{-i}\cdot), 1\leq i\leq r$. Then
$$
\widehat{\mu}_i(\xi)=\prod_{k=0}^\infty  M_N\big (u^k\rho^i\xi\big )
$$
for $1\le i\le r$. Then $\mu_\rho$ is the convolution of $\mu_i$,
$i=1, 2,\ldots, r$. Let $\Lambda$ be a bi-zero set of $\mu_{\rho,
N}$. We claim that $\Lambda$ is also a bi-zero set of $ \mu_i$ for
some $i$. Indeed, let $\lambda_j=\rho^{-k_jr-i_j}a_j/N$,  $1\le i_1,
i_2\le r, j = 1,2$, be any two distinct elements in $\Lambda$. The
bi-zero property of $\Lambda$ for $ \mu$ implies that
$$
\rho^{-k_1r-i_1}a_1/N - \rho^{-k_2r-i_2}a_2/N  = \rho^{-kr-i}a/N
$$
Without loss of generality assume $k_1, k_2 \geq k$, then we have
$ u^{(k_1-k)}\rho^{i_1-i} a_1 - u^{(k_2-k)}\rho^{i_2-i} a_2 =a$. This implies $i_1=i_2=i$ because the minimal polynomial of $\rho$ is $x^r-u$. Hence $\Lambda$ is a bi-zero set of $\mu_i$, and by Lemma
\ref{th2.2}, $\Lambda$ cannot be a spectrum of $\mu$.  \eproof

  \bigskip

Next we consider $\rho\not \in {\Bbb Q}^{1/r}, r>1$.  We need two
lemmas.

\medskip

\begin{lem} \label{lemm2.4} Suppose $\Lambda$ is an infinity bi-zero set of $\mu_{\rho, N}$ with $0\in\Lambda$. Then  $\rho \not \in {\Bbb Q}^{1/r}$
 for all $ r \geq 1$
implies that $\rho$ is an algebraic integer.
\end{lem}

\medskip

\proof Since $\Lambda\setminus\{0\}\subset {\mathcal
Z}(\widehat{\mu}_{\rho, N})$, we denote
$\Lambda=\{\lambda_k\}_{k=0}^\infty$ so that $\lambda_0=0$ and
$\lambda_k=\rho^{-n_k}\frac {a_k}N$, where $N\nmid a_k$ for $k\ge
1$. We can assume that  $n_k\le n_{k+1}$ for $k\ge 1$. Fix $\ell\ge
1$. For any integer $G>0$ and $k>\ell$, by the bi-zero property of
$\Lambda$, we have
$$
\lambda_k-\lambda_{\ell}=\rho^{-n_{k, \ell}}\frac{a_{k,  \ell}}N,
\quad a_{k,  \ell}\in \Z\setminus N\Z.
$$
We claim  $\#\{k:n_{k, \ell}\le G\} \leq (N-1)G$. Otherwise, by the
pigeon hole principle, there exist $k_1, k_2$ such that $n_{k_1,
\ell}=n_{k_2, \ell}\le G$ and $N\mid (a_{k_1, \ell}-a_{k_2, \ell})$.
Then, by the definition of ${\mathcal Z}(\widehat{\mu}_{\rho, N})$
and $\rho \not \in {\Bbb Q}^{1/r}$ for all $ r \geq 1$, we have
$$\lambda_{k_1}-\lambda_{k_2}=\lambda_{k_1}-\lambda_\ell+\lambda_\ell-\lambda_{k_2}=\rho^{-n_{k_1,
\ell}}\frac {a_{k_1, \ell}-a_{k_2, \ell}}N\not\in {\mathcal
Z}(\widehat{\mu}_{\rho, N}).$$
 Hence the claim follows. Taking any $k>\ell$ such that $n_{k, \ell}>n_\ell$, we conclude from
$$
\rho^{-n_k}\frac{a_k}N-\rho^{-n_\ell}\frac{a_\ell}N-=\rho^{-n_{k,
\ell}}\frac{a_{k, \ell}}N
$$
that there exists a polynomial $p(x)=a_\ell x^s+bx^t+c$ with $s>t$
and $p(\rho)=0$. Let $\varphi (x)=
c_0x^m+c_{1}x^{m-1}+\cdots+c_m\in\Z[x]$ be the minimal polynomial of
$\rho$. This implies that $\varphi(x)\mid p(x)$,  and thus $c_0\mid
a_\ell$. Let $\ell$ run through all $\lambda_\ell \in \Lambda$. Then
\begin{equation}\label{5.01}
\frac 1{c_0}\Lambda \setminus\{0\} \subseteq {\mathcal
Z}(\widehat{\mu}_{\rho, N}).
\end{equation}
To show that $\frac 1{c_0}\Lambda$ is a bi-zero set of $\mu_{\rho,
N}$ we need to prove that
\begin{equation}\label{5.02}
\frac 1{c_0}(\Lambda-\Lambda) \setminus\{0\} \subseteq {\mathcal
Z}(\widehat{\mu}_{\rho, N}).
\end{equation}
For any  $\lambda_{k_1}\ne\lambda_{k_2} \in \Lambda$,  by the claim
there exists $k$ such that $\min \{ n_{k,k_1}, n_{k,k_2} \} >
n_{k_1,k_2}$, thus
$$\rho^{-n_{k_1,k_2}}\frac {a_{k_1,k_2}}N
=\lambda_{k_1}-\lambda_{k_2}=(\lambda_{k_1}-\lambda_k)-(\lambda_{k_2}-
\lambda_k)=\rho^{-n_{k,k_2}}\frac
{a_{k,k_2}}N-\rho^{-n_{k,k_1}}\frac {a_{k,k_1}}N.$$ Similar to the
above, we have $c_0\mid a_{k_1,k_2}$. Then \eqref{5.02} holds.

By repeating the same argument, we see that $\frac 1{c_0^k}\Lambda$
is also a  bi-zero set of $\mu_{\rho, N}$ for  any $k\ge 1$. This
force $c_0=1$.\eproof
\bigskip

For any $x\in \R$, let $\|x\|=|\langle x\rangle|$,  where $\langle x\rangle$ is the unique number such that $\langle x \rangle \in (-1/2, 1/2]$ and $x-\langle x\rangle \in\Z$. Clearly $||x||$ is the distance from $x$ to ${\Bbb Z}$.

\medskip

\begin{lem}\label{lemm2.5} Let $\rho$ be a root of  $x^m+c_1x^{m-1}+\cdots+c_m\in\Z[x]$. Then for any $a\in\Z\setminus N\Z$,
\begin{equation} \label{new5.1}
\max_{1\le n\le m} \|\rho^{-n}\frac aN \|\ge\left(N\sum_{n=1}^m|c_n|\right)^{-1}:=\alpha>0 .
\end{equation}
\end{lem}

\medskip

\noindent {\bf Proof.}
Denote $\rho^{-n}\frac aN =\langle\rho^{-n}\frac aN\rangle + k_n$, \  $1\le n\le m$. Then
\begin{equation} \label{new5.2}
 \frac aN +\sum_{n=1}^m c_{n} \langle\rho^{-n}\frac aN\rangle + \sum_{n=1}^m c_n k_n=0.
\end{equation}
If $|\langle\rho^{-n}\frac aN\rangle|<\alpha$ for $1\le n\le m$, then $|\sum_{n=1}^m c_{n} \langle\rho^{-n}\frac aN\rangle|<\frac 1N$. This contradicts \eqref{new5.2} as $a\in\Z\setminus N\Z$. Hence the result follows.
\eproof

\bigskip

\begin{prop}\label{pro2.6}Let $\rho$ be an irrational and $\rho\not\in\Q^{1/r}$ for any $r>1$. Then $\mu_{\rho, N}$ is not a spectral measure.
\end{prop}

\medskip

\noindent {\bf Proof.} Suppose on the contrary that $\mu_{\rho, N}$
is a spectral measure. Then, by Lemma \ref{lemm2.4}, $\rho$ is an
algebra integer, and $\phi(x)=x^m+c_{1}x^{m-1}+\cdots+c_m\in\Z[x]$ is the minimal polynomial of $\rho$.

Let $\Lambda$ be a spectrum of $\mu_{\rho, N}$ with $0\in\Lambda$.
Denote $\Lambda_k=\Lambda\cap\{\rho^{-k}\frac aN: a\in\Z\setminus
N\Z\}$ for $k\ge 1$. Then $\#\Lambda_k\le N-1$ for $k\ge 1$ (by the
proof of Lemma \ref{lemm2.4}). Let $M_N(\xi)$ be the mask polynomial
and let  $G(\xi)=\sum_{i=1}^{N-1}|M_N(\xi+\frac iN)|^2$. Then by
applying Theorem \ref{th2.1}  to the point mass measure $\frac 1N
\delta_{\{ 0, \cdots , N-1\}}$, we have
$$
G(\xi)+|M_N(\xi)|^2=\sum_{i=0}^{N-1}|M_N(\xi+\frac iN)|^2 =1,
$$
and
 hence $G(0)=0$. Observing that $G(z)$ is an entire function, then there exists an entire function $H(z)$ and integer $t>0$ such that $G(z)=z^tH(z)$ and $H(0)\ne 0$. To prove that $Q(\xi)=|\widehat{\mu}_{\rho, N}(\xi)|^2+\sum_{k=1}^\infty
\sum_{\lambda\in\Lambda_k}|\widehat{\mu}_{\rho,
N}(\xi+\lambda)|^2\not\equiv 1$, we first observe that  for any
$\xi$,
\begin{eqnarray}
\sum_{\lambda\in\Lambda_k}|\widehat{\mu}_{\rho, N}(\xi+\lambda)|^2
&=&\sum_{\lambda\in\Lambda_k} \prod_{j=1}^k |M_N(\rho^j(\xi+\lambda))|^2 \cdot |\widehat{\mu}_{\rho, N}
(\rho^k(\xi+\lambda))|^2 \label{ine5.2}\\
&\le& \sum_{\lambda\in\Lambda_k}|M_N(\rho^k(\xi+\lambda))|^2 \nonumber\\
&\leq & \ G(\rho^k\xi).\nonumber
\end{eqnarray}
 (The last inequality follow from $ \lambda\in\Lambda_k$,\
$\rho^k\lambda=\frac aN \not = 0$ ,  $a\not \mid N$).  Let $m$ and
$\alpha (< 1/2)$ be defined as in Lemma \ref{lemm2.5}, and let
 $\beta=\min \{1-|M_N(x)|^2:{\alpha/2\le |x|\le 1-\alpha/2}\}$. Then obviously $\beta>0$. Note that for each  $k>m$ and $\lambda\in \Lambda_k$,
 $$
 \rho^{j}\lambda=\rho^{-(k-j)}\frac aN, \quad j=1,2,\ldots, k-1.
 $$
  Hence for $0\le \xi\le \alpha/2, k>m$, by Lemma \ref{lemm2.5}, there exist $k-m\leq \ell_\lambda \leq k-1 $ such that
 $||\rho^{\ell_\lambda}(\xi +\lambda)||^2\geq \alpha/2$.
 Hence from (\ref {ine5.2}), we have
\begin{eqnarray} \label{ine5.3}
\sum_{\lambda\in\Lambda_k}|\widehat{\mu}_{\rho,N}(\xi+\lambda)|^2
&\le & \sum_{\lambda\in\Lambda_k}
|M_N(\rho^{\ell_\lambda}(\xi+\lambda))|^2 \cdot |M_N
(\rho^k(\xi+\lambda))|^2  \nonumber\\
&\le &  (1-\beta) \sum_{\lambda\in\Lambda_k}   |M_N(\rho^{k}(\xi+\lambda))|^2  \nonumber \\
&\le& (1-\beta) G(\rho^k\xi) \nonumber.
\end{eqnarray}

\medskip
 Note that $\Lambda\setminus\{0\}=\cup_{k\in\N}\Lambda_k$, and $\Lambda_{k_1}\cap\Lambda_{k_2}=\emptyset$ when $k_1\ne k_2$ since
 $\lambda \not\in \Q^{\frac 1r} $ for all $r\in\N$. Hence, by \eqref{ine5.2} and  \eqref{ine5.3},
\begin{eqnarray} \nonumber
Q(\xi)
& = & \sum_{\lambda\in\Lambda}|\widehat{\mu}_{\rho,N}(\xi+\lambda)|^2 \\ \label{new5.6}
& = &|\widehat{\mu}_{\rho,N}(\xi)|^2+\sum_{k=1}^\infty\sum_{\lambda\in\Lambda_k}|
\widehat{\mu}_{\rho, N} (\xi+\lambda)|^2\\ \nonumber
&\le& |\widehat{\mu}_{\rho,N}(\xi)|^2+\sum_{k=1}^{m}G(\rho^k\xi)+(1-\beta)
\sum_{k>m} G(\rho^k\xi).
\end{eqnarray}

 On the other hand, recall that $G(z)=z^tH(z)$ and $H(0)=0$, then $0<C_1\le |H(z)|\le C_2$ if $|z|\le
 \eta \le \alpha/2$ for some small $\eta$. Therefore for $0\le\xi\le \eta$,
\begin{eqnarray}\label{add}
\frac{C_1\rho^{m t}}{1-\rho^t}\xi^t\le\sum_{k=m}^\infty G(\rho^k\xi)\le\frac{C_2\rho^{mt}}{1-\rho^t}\xi^t,
\end{eqnarray}
and
\begin{eqnarray} \nonumber
 |\widehat{\mu}_{\rho, N}(\xi)|^2
 & = & \prod_{k=1}^\infty |M_N(\rho^k\xi)|^2=\prod_{k=1}^\infty\big(1-G(\rho^k\xi)\big) \\ \label{new5.7}
 &\le&
 e^{-\sum_{k=1}^\infty G(\rho^k\xi)}\le 1-\sum_{k=1}^\infty G(\rho^k\xi)+o\left(\sum_{k=1}^\infty G(\rho^k\xi)\right),
\end{eqnarray}
 where $o(\xi)$ satisfies that $\lim_{\xi\rightarrow 0}
 o(\xi)/\xi=0$.
Hence, by \eqref{new5.6} and \eqref{new5.7}, we have
\begin{equation} \label{ine5.5}
Q(\xi)\le 1-\beta\sum_{k=m}^\infty G(\rho^k\xi)+o\left(\sum_{k=1}^\infty G(\rho^k\xi)\right).
\end{equation}
  By \eqref{add} this  implies $Q(\xi)<1$ for $\xi>0$
  small enough. That $\Lambda$ cannot be a spectrum follows by Theorem \ref{th2.1}.
\eproof

\bigskip

In view of Propositions \ref{pro2.3} and \ref{pro2.6}, we have to prove that  $\mu_{\rho, N}$  cannot be a spectral measure in the  remaining cases (iii) and (iv) in \S 1 for $\rho =p/q $. These will be proved in the remaining sections.

\bigskip
\bigskip

\section {\bf Structure of bi-zero sets for rational contraction}

In this section we will consider $\rho =  p/q$, we assume $p, q$ are co-primes throughout. Let $\Lambda=\{\lambda_k\}_{k=0}^\infty\subseteq {\mathcal Z}(\widehat{\mu}_{\rho,N})$ (with $\lambda_0=0)$ be a bi-zero set of
$\mu_{\rho, N} $  Then by
(\ref{eq2.1}),
\begin{equation} \label {eq2.2}
\lambda_k=\Big (\frac qp \Big )^{n_k}\frac {a_k}{N} \quad \hbox {with} \ \
a_k \in \Z \setminus N\Z, \quad k\geq 1.
\end{equation}
In the following, we will give another expression of the $\lambda_k$ which is more
convenient to use here.

\medskip

\begin {lem} \label {lemm3.1}
Let $\Lambda$ be a bi-zero set of $\mu_{\rho, N}$ with $\rho=\frac
pq$. Then  there exists $m_0>0$ such that  each $\lambda_k \in \Lambda \setminus \{0\}$ admits an
expression
\begin{equation}\label{eq2.3}
 \lambda_k=p^{-m_0}q^{m_k}\frac {c_k}N\quad \hbox {with} \ \  c_k \in \Z \setminus q\Z \ \  \hbox {and} \ \ m_k \geq m_0
 \end{equation}
(note that $N$ can be a factor of $c_k$). Moreover, if $N\mid q$, then we can write
\begin{equation*}
 \lambda_k=p^{-m_0}q^{m_k}\frac {c_k}N\quad \hbox {with} \ \   c_k \in \Z \setminus N\Z \ \  \hbox {and} \ \ m_k \geq
 m_0.
 \end{equation*}
\end{lem}

\medskip
\noindent {\bf Proof.} For the expression of $\lambda_k$ in
(\ref{eq2.2}), we let $a_k=a_k'q^{l_k}$ so that $q\nmid a_k'$. Then
we can write $\lambda_k$ as
\begin{equation}\label{eq}
\lambda_k=\Big(\frac qp \Big)^{n_k+l_k}\frac{a_k'p^{l_k}}N:=\Big(\frac
qp\Big)^{m_k}\frac {b_k}N \ ,
\end{equation}
where $q$ is not a factor of $b_k$ for $k\ge 1$.  Let $m_0\geq 1$ be the smallest among all such $m_k$, and denote the corresponding
$\lambda_i\in \Lambda$ by $(\frac qp)^{m_0}\frac {b_i}N $. Then by the bi-zero property, for any $m_k>m_0$,
\begin{eqnarray*}\label{equ1}
\Big(\frac qp\Big)^{m_0}\frac {b_i}N\ -\ \Big(\frac
qp\Big)^{m_k}\frac {b_k}N\ =\ \Big(\frac qp\Big)^{m}\frac {b}N.
\end{eqnarray*}
It is easy to see that $m=m_0$, and then $p^{m_k-m_0}$ is a factor of $b_k$. It follows from this that we can rewrite $\lambda_k$ as
\begin{equation*}
 \lambda_k=p^{-m_0}q^{m_k}\frac {c_k}N,
 \end{equation*}
  where $q\not\mid c_k$ for $k\ge 1$.

  \medskip

  The second assertion follows by observing that the $l_k$ in
  \eqref{eq} is zero  (as  $q\nmid a_k$ follows by $N\mid q$ and $N \nmid a_k$). Hence the above $c_k = a_k/p^{m_k -m_0}$ is not divisible by $N$ by  \eqref{eq2.2}.\eproof

\bigskip

\begin{coro} \label {th2.4}
Let $\Lambda$ be a bi-zero set of $\mu_{\rho, N}$ and let $N\mid q$
and $\rho=\frac pq$. Denote $Q=\{q^ma : \  a\in {\Bbb Z}\setminus N{\Bbb Z}, m\geq 0\}$.
Then
\begin{eqnarray}\label{eq2.4}
(\Lambda-\Lambda)\setminus\{0\}\subseteq\frac 1{\rho^{m_0}N}Q
\subset {\mathcal Z}(\widehat{\mu}_{\rho, N}(\xi)),
\end{eqnarray}
where $m_0$ is as in Lemma \ref{lemm3.1}.
\end{coro}

\medskip

\noindent {\bf Proof.} It suffices to show that
\begin{eqnarray*}
(\Lambda-\Lambda)\setminus\{0\}\subseteq  \frac 1{p^{m_0}N} \
\Big\{q^ma: m\ge m_0, \  a\in\Z\setminus N\Z\Big \}.
\end{eqnarray*}
 If
$m_k> m_l$ for $k\ne l$, we have
$\lambda_k-\lambda_l=p^{-m_0}q^{m_l}\frac {q^{m_k-m_l}c_k-c_l}N\in
\frac 1{\rho^{m_0}N}Q$ because $N\nmid c_l$. If $m_k=m_l$ for $k\ne
l$, then by Lemma \ref{lemm3.1},
$$
\lambda_k-\lambda_l=p^{-m_0}q^{m_k}(c_k-c_\ell)/N=p^{-m_0}q^{m_k+\alpha}c/N
$$
where $q\nmid c$. By the bi-zero property in (\ref{eq2.2}), we have
$$
p^{-m_0}q^{m_k+\alpha}c/N=\lambda_k-\lambda_l=(\frac qp)^na/N
$$
where $N \nmid a$. Then $q^{m_k+\alpha-n}c=p^{m_0-n}a$, which
implies that $m_k+\alpha=n$, and thus $a=cp^{\alpha+m_k-m_0}$. Hence
$N\nmid c$ and the claim follows.  \qquad $\Box$

\bigskip

It is well-known that every positive integer has a unique $q$-adic expansion. In order to do this for  all integers in $\Z$, we use the $q$-adic expansion on
the set $ \C=\{-1, 0, \cdots q-2\}$. In the following, we will establish a relation of
the $\lambda_k$ in the bi-zero set $\Lambda$ with such expansion. We characterize the maximal bi-zero set by certain tree-structure.  We need the addition condition that $N\mid q$, and a special selection map to be defined in the following.

\bigskip

 Let $\Omega_N =\{0, \cdots N-1\}$ and let
$\Omega^*_N=\bigcup_{k=0}^\infty \Omega^k_N$ be the set of finite
words (by convention $\Omega^0_N=\{\emptyset\}$). We use ${\bf i} =
i_1 \cdots i_k$ to denote an element in $\Omega_N^k$, and $|{\bf
i}|=k$ is the length. For any ${\bf i}, {\bf j} \in\Omega^*_N$,\
${\bf i}{\bf j}$ is their natural conjunction. In particular,
$\emptyset{\bf i} ={\bf i}$, ${\bf i} 0^\infty={\bf i} 0 0\cdots$
and $0^k=0\cdots 0\in \Omega_N^k$.

\bigskip

\begin{defi} \label {def4.1}{\it  Suppose $N,q$ are positive integers and $N\mid q$.
 We call  a map $\iota: \Omega_N^*\rightarrow \{-1,0,...,q-2\}$ a {\it selection mapping} if

\vspace{0.1cm}

\ \ {\rm (i)} $\iota(\emptyset) =\iota(0^n) =0$  for all $n\geq1$;

\vspace {0.1cm}

\ {\rm (ii)} for any ${\bf i} =i_1\cdots i_k \in \Omega_N^k$,\ \
$\iota({\bf i}) \in (i_k + N\Z)\cap\C$, where $\C=\{-1, 0, 1,
\ldots, q-2\}$;

\vspace{0.1cm}

{\rm (iii)}  for any  ${\bf i} \in \Omega_N^{\ast}$, there exists
${\bf j}\in \Omega_N^{\ast}$ such that  $\iota$ vanishes eventually
on ${\bf ij}0^{\infty}$, i.e., $\iota( {\bf i}{\bf j}0^k)=0$ for
sufficient large $k$.}
\end{defi}

\medskip

Note that $\C \equiv \Omega_N \oplus N\{0, \cdots, r-1\} (\hbox
{mod} \ q)$ where $q=rN$, and $\iota$ is a selection map on each
level $k$.  More explicitly, (ii) means
\begin{equation} \label {eq4.1}
\iota({\bf i}) =\left\{
\begin{array}{ll}
i_k + Nt, & \hbox{if} \ \ 0\leq i_k \leq N-2, ; \\
i_k +Nt', & \hbox{if} \ \  i_k = N-1,
 \end{array}
 \right.
\end{equation}
where $t\in \{0, \cdots r-1\}$ and $t' \in \{-1, 0, \cdots, r-2\}$.

\bigskip

Next we let
$$
\Omega^{\iota}_N=\{\ {\bf i} =i_1\cdots i_k \in \Omega_N^*: i_k\ne
0,
 \ \  \iota({\bf i} 0^n)=0 \,\,\mbox{for
sufficient large $n$}\} \cup \{\emptyset\}
$$
and for any ${\bf i}\in\Omega^{\iota}_N$ we define
$$
\iota^*({\bf i})=\sum_{n=1}^\infty\ \iota({\bf i} 0^\infty)|_n)\ q^{n-1},
$$
Here  we regard  ${\bf i}0^\infty = {\bf i} 00 \cdots $, and ${\bf
i}0^\infty|_n$ denotes the word of the first $n$ entries. Clearly
$\iota^*(\emptyset)=0$.

\bigskip

Let
 $Q =\{q^ma :\  a \in \Z\setminus N\Z, m \geq 0 \}
 $
 be as in Corollary \ref{th2.4}, a subset $ L\setminus\{0\} \subset Q$ is called a {\it D-set} of $Q$ if  $0\in L$ and $L-L \subset Q\cup
 \{0\}$  ($D$ for difference),   and call it a {\it maximal} D-set if for any
$n\in Q \setminus L$,  $L \cup \{n\}$ is not a D-set.  The main idea
of the proof of the following theorem is in \cite {DHL} (and the selection map is called a maximal map
there)(see also \cite{DHS}). We provide a simplified proof here for completeness.

\bigskip

\begin{prop} \label{th4.2} Suppose $N\mid q$. Then
$L \subset Q:=\{q^ma: m\ge 0, a\in\Z\setminus N\Z\}$ is a maximal
D-set of $Q$ if and only if $L =\iota^*( \Omega^{\iota}_N)$ for some
selection map $\iota$.
\end{prop}

\medskip

{\noindent \bf  Proof.} We first prove the sufficiency. For a
selection map $\iota$,  it is direct to check that $L
=\iota^*(\Sigma^{\iota}_N)$ is a D-set of $Q\subseteq \Z$ by the
definition of $\iota$. We need only show that $L$ is maximal in $Q$.
Suppose otherwise, there exists $n\not\in L$ and $L\cup\{n\}$ is a
D-set. We can express $n$ uniquely as
\begin{equation} \label{eq4.2}
n=a_0+a_1q+\cdots+a_\ell q^\ell, \quad  a_i\in \C=\{-1, 0, 1,
\ldots,
 q-2\}.
\end{equation}
 We claim that $a_0=\iota(i_1)$ for some $i_1 \in \Omega_N$. If otherwise, let $j \in \Omega_N  =\{0, \cdots , N-1\}$ such that $a_0 \in
j + N \Z$. In view of property (ii) of $\iota$ (or (\ref{eq4.1})),
$N\mid(a_0-\iota(j))$. By property (iii) of $ \iota$, there exists
${\bf i}=i_1\cdots i_k\in\Omega^{\iota}_N$ with $i_1=j$. Then
$$
n-\iota^*({\bf i})=a_0-\iota(j)+qb,
$$
where $b$ is an integer. Hence $n-\iota^*({\bf i})  \not\in Q$ (as
it has a factor $N$, and not a factor of $q$). This contradicts
that $L \cup\{n\}$ is a D-set of $Q$, and the claim follows.

\vspace{0.1cm} Similarly, by considering $n- \iota(i_1) =n- a_0$ in
(\ref{eq4.2}),  we can show that $a_1=\iota(i_1i_2)$ for some $0\le
i_2<N-1$, and so on. After finitely many steps, we have
$n=\iota^*({\bf i})$ for some ${\bf i}\in\Omega^{\iota}_N$, which
contradicts $n \not \in L$, and the sufficiency  follows.

\medskip

Conversely, suppose that $L$ is a maximal D-set of $Q$. Denote
$L=\{\lambda_k\}_{k=0}^\infty$ with $\lambda_0=0$. Then $\lambda_k$
can be expressed by
$$
\lambda_k =a_{k,0}+a_{k,1}q + \cdots+a_{k, l_k}q^{l_k}
=\sum_{n=0}^\infty a_{k, n}q^{n},
$$
where $-1\le a_{k, n}\le q-2$ for $0\le n\le l_k$ and
$a_{k, n}=0$ for $n>l_k$. Note that all
$a_{0, n}$ are zero. We first consider $ \{a_{k,0}: k \geq 0\}$, the first coefficients of the $\lambda_k$'s.
As $a_{k, 0}$ can be written uniquely  as $i_k +
N\alpha_k\in\C=\{-1, 0, \ldots, q-2\}$ for some $i_k \in \Omega_N =
\{0, \cdots, N-1\}$. We claim that
\begin{eqnarray}\label{theta}
\{a_{k, 0}: k \geq 0\} =\{i+N\alpha_i: i\in \Omega_N \}\
\subseteq \C.
\end{eqnarray}
(Here $\alpha_i$ depends only on $i$, but not on $k$, hence the set
have $N$ elements.) Indeed if $\{a_{k, 0}: k \geq 0\}\supsetneqq\{i+N\alpha_i: i\in \Omega_N \}$, then
there exist $k_1$ and $k_2$ such that $N\mid(a_{k_1, \, 0}-a_{k_2,
\, 0})$. Hence
$$
\lambda_{k_1}-\lambda_{k_2}=a_{k_1, \, 0}-a_{k_2, \, 0}+qb\not\in Q,
$$
(same reasoning as the above) which contradicts that $L$ is a  D-set
in $Q$. If $\{a_{k, 0}: k \geq 0\}\subsetneqq\{i+N\alpha_i:
i\in \Omega_N \}$,  then there exists $0\le i'\le N-1$ such that
$N\nmid (a_{k, \, 0}-i')$ for $k\ge 0$. Clearly $L\cup\{i'\}$ is a
D-set in $Q$, which contradicts the maximality of $L$. This
proves the claim. We rewrite \eqref{theta} as
$$
\{a_{k, 0}: k \geq 0\} =\{i_0+N\alpha_{i_0,0}: i_0\in \Omega_N \}\
\subseteq \C.
$$

\medskip

From the claim, we can define $\iota$ on $\Omega_N$ by $
\iota(i)=i+N\alpha_{i, 0}, \ i=0, 1, \ldots, N-1$ and in particular
$\iota(0)=0$.  Similarly we can show that, for each $0\le i_0\le
N-1$,
$$
\{a_{k,1}:\ a_{k,0}=i_0+N\alpha_{i_0,
0}\}_{k=0}^\infty=\{i_1+N\alpha_{i_1, 1}: i_1\in \Omega_N\}\subseteq
\C
$$
and define $\iota(i_0i)=i+N\alpha_{i, 1}, \ i=0, 1, \ldots, N-1.$
Again, we can show that
$$
\{a_{k,2}:\ a_{k,0}=i_0+N\alpha_{i_0, 0}\,\, \mbox{and}\,\,
a_{k,1}=i_1+N\alpha_{i_1, 1}\}_{k=0}^\infty=\{i_2+N\alpha_{i_2, 2}:
i_2\in \Omega_N\}
$$
and define $\iota(i_0i_1i)=i+N\alpha_{i, 2}, \ i=0, 1, \ldots, N-1.$
Inductively, we can define a map $\iota$ on $\Omega_N^*$ (with
$\iota(\emptyset)=0$).  By the construction of $\iota$, it is easy
to see that (i) and (ii) in Definition \ref{def4.1} are satisfied.
For any ${\bf i}=i_0i_1\cdots i_n\in \Omega_N^*$ with $i_n\ne 0$,
again by the construction of $\iota$, there exists infinitely many
$\lambda_k$ such that $a_{k, t}=i_t+N\alpha_{i_t, t}$ for $0\le t\le
n$. Fix such a $k$, if $k\ge l_k$, we have $ \lambda_k
=\sum_{n=0}^\infty a_{k, n}q^{n}=\iota^*({\bf i}) $; If $k<l_k$,
 there exists ${\bf j}=j_{n+1}j_{n+2}\cdots j_{l_k}$
such that $a_{k, t}=\iota(i_0\cdots i_n j_{n+1}\cdots
j_t)$ for $n+1\le t\le l_k$. Then
$$
\lambda_k =\sum_{n=0}^\infty a_{k, n}q^{n}=\iota^*({\bf ij}).
$$
This implies that (iii) in Definition \ref {def4.1} holds. Hence,
$\iota$ is a selection mapping and $L\subseteq
\iota^*(\Omega_N^\iota)$. The necessity follows by the maximal
property of $L$ and the  proof of the sufficiency. \eproof

\bigskip
It follows directly from Corollary \ref {th2.4} and Proposition
\ref{th4.2} that

\medskip

\begin{theo} \label {thm3.5} Suppose $\rho =p/q$ and $N\mid q$. Then $\Lambda \subset {\mathcal Z}(\widehat{\mu}_{\rho, N})$
is a maximal bi-zero set if and only if there exist $m_0\ge 1$ and a
selection map $\iota$ such that
$\Lambda=\rho^{-m_0}N^{-1}(\iota^*(\Omega^{\iota}_N))$.
\end{theo}

\bigskip
In particular, we see that for $p=1$, the spectrum $\Lambda$ in
(\ref{eq3.1}) corresponding to the case $m_0=1$  and the selection
map $\iota$ is to take $\iota ({\bf i}) =i_k$ in (\ref{eq4.1}).  Also we observe that $\iota^*(\Omega^{\iota}_N)$ is an infinite set, we have

\medskip
\begin{coro} Suppose $\rho =p/q$ and $N\mid q$, then $L^2(\mu_{\rho, N})$ admits an infinite exponential orthonormal set.
\end{coro}

\bigskip

\section{\bf Spectrality for rational contraction}

\setcounter{equation}{0}

In this section, we prove the necessity of Theorem \ref{th1.1}  when  $\rho$ is a rational number.

\bigskip

\begin{prop}\label{prop3.1} Let $\rho=\frac pq$ and $1\le\gcd(N, q)<N$, then $\mu_{\rho, N}$ is not a spectral
measure.
\end{prop}

\proof Suppose on the contrary that $\mu_{\rho, N}$ is a spectral measure. Let $\Lambda$ be a spectrum of $\mu_{\rho, N}$ with $0\in \Lambda$.
Denote $d=\gcd(N, q)$. If $d=1$,  by Lemma \ref{lemm3.1}, we have
$$\Lambda\subseteq p^{-m_0}\{q^m\frac aN: m\ge m_0, a\in\Z\setminus q\Z\}\cup\{0\}.$$
Denote $\D=\{0, 1, \ldots, N-1\}$ and let
$\mu'=\delta_{\rho^1\D}\ast \delta_{\rho^2\D}\ast\cdots \ast
\delta_{\rho^{m_0}\D}\ast \delta_{\rho^{m_0+2}\D}\ast\cdots$ be the
convolution of $\delta_{\rho^k\D}$ for $k\ge 1$ and $k\ne m_0+1$
(here $\delta_A=\frac 1{\#A}\sum_{a\in A}\delta_a$ and $\delta_a$ is
the Dirac measure). Then $\mu_{\rho,
N}=\delta_{\rho^{m_0+1}\D}\ast\mu'$. We claim that $\Lambda$ is a
bi-zero set of $\mu'$. The claim leads to a contradiction by Lemma
\ref{th2.2}. We prove the claim by assuming that $\rho^{-m_0-1}\frac
aN\in \Lambda-\Lambda$ where $a\in\Z\setminus N\Z$. Then there exist
$k, l$ such that
 $$\rho^{-m_0-1}\frac aN=p^{-m_0}q^{-m_k}\frac {a_k}N-p^{-m_0}q^{-m_l}\frac {a_l}N,$$
 where $a_k, a_l\in (\Z\setminus q\Z)\cup\{0\}$. Then $p\mid a$. Hence $\rho^{-m_0-1}\frac aN=\rho^{-m_0}\frac{qa/p}N\in{\mathcal Z}(M_{\rho^{m_0}N})$ and the claim follows;  If $1<d<N$,  write $N=N'd, q=q'd$. Then $\D=\C+d\E$, where $\C=\{0,
1, \ldots, d-1\}$ and $\E=\{0, 1, \ldots, N'-1\}$. Note that\
$M_N(\xi)=M_d(\xi)M_{N'}(d\xi)$\ and
$$
\widehat{\mu}_{\rho, \, N}(\xi)=\prod_{k=1}^\infty
M_N(\rho^k\xi)=\prod_{k=1}^\infty M_d(\rho^k\xi)\prod_{k=1}^\infty
M_{N'}(\rho^kd\xi).
$$
Let $\nu$ be the probability measure such that
$$
\widehat{\nu}(\xi)=\prod_{k=1}^\infty M_d(\rho^k\xi)\prod_{k\ge 1,
\, k\ne m_0+1} M_{N'}(\rho^kd\xi).
$$
Then $\mu=\nu\ast \delta_{\rho^{m_0+1}d\E}$ . We
claim that $\Lambda$ is a bi-zero set of $\nu$. Hence the proposition
follows by Lemma \ref{th2.2} again.

\medskip

To prove the claim, we let $\eta\in (\Lambda-\Lambda)
\setminus\{0\}\ (\subset {\mathcal Z}(\mu_{\rho, N})$), then either
$\eta \in {\mathcal Z}(\nu)$ or $\eta \in {\mathcal
Z}\big(M_{N'}(\rho^{m_0+1} d(\cdot))\big)$. The first case satisfies the
claim trivially. Hence we need only consider the second case, i.e.,
there exists $\eta\in (\Lambda-\Lambda)$ such that $\eta\in
{\mathcal
 Z}\big(M_{N'}(\rho^{m_0+1} d(\cdot))\big)$. By (\ref{eq2.1}),
 we have
 $\eta =\frac 1{\rho^{m_0+1}d}\frac a{N'}$ $\left(=(\frac qp)^{m_0+1}\frac aN \right)$ with $N'\nmid a$;  also by  \eqref{eq2.3},
 there exist $k, \ell$ such that
 $$
 \eta=p^{-m_0}q^{m_k}\frac {c_k}N-p^{-m_0}q^{m_\ell}\frac {c_\ell}N,
 $$
 where $q\nmid c_k$ and $q\nmid c_\ell$. Hence we have
 $\frac {q^{m_0+1}a}{p}=q^{m_k} {c_k}-q^{m_\ell} {c_\ell}$.
This implies that $p\mid a$. By letting $a'=qa/p$, we see that
$N\nmid a'$ (as $N'\nmid a$).  Therefore
$$
\eta=(\frac qp)^{m_0+1}\frac a{N}=(\frac qp)^{m_0}\frac {a'}{N} \in
{\mathcal Z}\big(M_N(\rho^{m_0} (\cdot))\big).
$$
As ${\mathcal Z}\big(M_N(\rho^{m_0} (\cdot))\big)\subset {\mathcal
Z}(\nu)$,  the claim follows. \eproof

\bigskip

\begin{prop}\label{prop3.2} Let $\rho=\frac pq$. If $q\mid N $ and $p> 1$, then $\mu_{\rho, N}$ is not a
spectral measure.
\end{prop}

\bigskip

\medskip

The proof of this case is more elaborate.  We show that any maximal
bi-zero set $\Lambda$ of $ \mu_{\rho, N}$ does not satisfy the
condition on $Q(\xi)$ in Theorem \ref{th2.1}(ii). To this end, we
define
$$
\mu_n=\delta_{\rho\Omega_N}\ast \cdots \ast
\delta_{\rho^{n}\Omega_N}
$$
for $n\ge 1$. Then
$$
\widehat{\mu}_n(\xi)=\prod_{i=1}^n
M_N(\rho^i\xi) \quad \hbox {and} \quad \widehat{\mu}_{\rho, \,
N}(\xi)=\widehat{\mu}_n(\xi)\widehat{\mu}_{\rho, \, N}(\rho^n\xi).
$$
We need a few technical lemmas.

\medskip

\begin{lem}\label{th4.4} Let $\iota$ be a selection mapping. Then
$$
\sum_{{\bf i} \in \Omega^{\iota}_N,\, |{\bf i}|\le
n}|\widehat{\mu}_{n+m_0-1}\big(\xi+\rho^{-m_0}N^{-1}\iota^*({\bf i}
)\big)|^2\le 1
$$
for $n\ge 1$ and $\xi\in\R$.
\end{lem}

\medskip

\noindent {\bf Proof.} First we prove the case for $m_0=1$.
According to the Bessel inequality, it suffices to show that
$\rho^{-1}N^{-1}\iota^*(\{{\bf i}\in\Sigma^{\iota}_N: |{\bf i}|\le
n\})$ is a bi-zero set of $\mu_n$. For any ${\bf i}, {\bf
 j} \in \Sigma^{\iota}_N$, ${\bf i} \ne{\bf j}$ and  $1\le |{\bf i}|,\  |{\bf j}|\le n$,  we let ${\bf i}'={\bf i} 0^{n-|{\bf i}|}:=i_1'\cdots
i_n'$ and ${\bf j}'={\bf j} 0^{n-|{\bf j}|}:=j_1'\cdots j_n'$. Let
$s$ be the smallest integer such that $i_s'\ne j_s'$. Then $s\le n$
and
$$
\iota^*({\bf i})-\iota^*({\bf j} )=\big(\iota({\bf i}'|_s) -\iota({\bf i}'|_s)\big)q^{s-1}+\alpha q^s
$$
for some integer  $\alpha$. By (\ref {eq4.1}), $\iota({\bf i}'|_s) -\iota({\bf i}'|_s)$ is not divisible by $N$. It follows from (\ref{eq2.0}) that,
$$
M_N\big (\rho^{s}\rho^{-1}N^{-1}\big(\iota^*({\bf i})-\iota^*({\bf
j})\big)\big ) =M_N\Big (\frac{p^{s-1}\big(\iota({\bf i}'|_s) -\iota({\bf i}'|_s)\big)}{N}\Big )=0\ .
$$
 This implies that $\widehat{\mu}_n(\rho^{-1}N^{-1}\big(\iota^*({\bf i})-\iota^*({\bf i})\big))=0$. Similarly, we have
$\widehat{\mu}_n(\rho^{-1}N^{-1}\iota^*({\bf i}))=0$ for any ${\bf
i}\in\Sigma^{\iota}_N$ and $0<|{\bf i}|\le n$. By Theorem
\ref{th2.1},
$$
\sum_{{\bf i} \in \Omega^{\iota}_N,\, |{\bf i}|\le
n}|\widehat{\mu}_{n}\big(\xi+\rho^{-1}N^{-1}\iota^*({\bf i}
)\big)|^2\le 1
$$
This completes the proof for $m_0=1$.  For $m_0>1$, we observe that
$$
|\widehat{\mu}_{n+m_0-1}(\xi)|=|\widehat{\mu}_{m_0-1}(\xi)|
|\widehat{\mu}_{n}(\rho^{m_0-1}\xi)|\le
|\widehat{\mu}_{n}(\rho^{m_0-1}\xi)|
$$
and apply the inequality. The result follows. \eproof

\bigskip

The following lemma is  a simple generalization of Lemma 2.10
in \cite{D}.
\begin{lem}\label{th4.5} Let $a=\ln p/\ln q$. Then for any  $\xi>1$ there exists $\xi'$ such that $\rho^2\xi^a \le \xi'\le \rho\xi$
and
$$
|\widehat{\mu}_{\rho, N} (\xi)|\ \le \ c\ |\widehat{\mu}_{\rho,
N}(\xi')|,$$
 where $c=\max\big \{\ |M_N(\xi)|:\  \frac 1{2q}\le
\xi\le \frac 12\big \} < 1$.
\end{lem}

\medskip

\noindent {\bf Proof.}  For any $x\in \R$, denote the unique number
$\langle x\rangle$ that satisfies $\langle x \rangle \in (-1/2,
1/2]$ and $x-\langle x\rangle \in\Z$. If $\langle \rho\xi \rangle
\not\in (-\frac 1{2q}, \frac 1{2q})$, then
$$
|\widehat{\mu}_{\rho, N}(\xi)|=|M_N(\rho\xi)|\ |\widehat{\mu}_{\rho,
N}(\rho\xi)|=|M_N(\langle \rho\xi)\rangle )|\ |\widehat{\mu}_{\rho,
N}(\rho\xi)|\le c|\widehat{\mu}_{\rho, N}(\rho\xi)|.
$$
Hence we obtain the desired inequality by letting $\xi'=\rho\xi $;
\medskip
If  $\langle\rho\xi\rangle \in (-\frac 1{2q}, \frac 1{2q})$, then
\begin{eqnarray}\label{eq4.4}
\rho\xi-\langle \rho\xi\rangle  =r_tq^t+\cdots+r_lq^l,
\end{eqnarray}
where  $ 0\le r_j<q$ for $t \leq j \leq l$ and $r_t>0$.  Then
$$
\langle \rho^{t+2}\xi\rangle =\big \langle \ \rho^{t+1}\langle
\rho\xi\rangle +\frac{r_{t}p^{t+1}}q \ \big \rangle.
$$
 Note that $|\rho^{t+1}\langle \rho\xi\rangle |<\frac 1{2q}$  and $\frac
1q\le|\langle \frac{r_{t}p^{t+1}}q\rangle | \le \frac{q-1}q$, then $
\langle \rho^{t+2}\xi\rangle \not\in (-\frac 1{2q}, \frac 1{2q}). $
\vspace {0.1cm}

By \eqref{eq4.4}, we have $\xi\ge q^{t},$
which implies  $\rho^{t}\ge \xi^{a-1}$, where $a=\ln p/\ln q$.
Let $\xi'=\rho^{t+2}\xi$, then $\xi'\ge \rho^2 \xi^a$, and hence
\begin{eqnarray*}
|\widehat{\mu}(\xi)|&=&|M_N(\rho\xi))|\ \cdots \
|M_N(\rho^{t+2}\xi)|\
|\widehat{\mu}(\rho^{t+2}\xi)|\\
&\le &|M_N(\langle\rho^{t+2}\xi\rangle)|\  |\widehat{\mu}(\rho^{t+2}\xi)|\\
&\le &c|\widehat{\mu}(\xi')|.
\end{eqnarray*}\eproof

\bigskip

\begin{lem}\label{th4.6}  Assume $p>1$, then there exist integers $b\geq 2$, $n_0\geq 2$, and real number  $\beta>1$, $C>1$ such that for any ${\bf i}
\in\Omega_N^{\iota}$ with ${n^b}<|{\bf i} |\le {(n+1)^b}$, $n \geq
n_0$, we have
$$
\Big|\widehat{\mu}_{\rho, N}\Big(\rho^{(n+1)^b+(m_0-1)}\big
(\xi+ \rho^{-m_0}N^{-1} {\iota^*({\bf i})}\big )\Big )
\Big |\le \frac C {n^\beta}
$$
for  $0\le \xi\le \frac1{2\rho^{m_0} N}$.
\end{lem}

\medskip

\noindent {\bf Proof.}  Note that $p>1$ implies that $q>2$. Let $b$ be an integer such that $b > 1+\frac
{\log a}{\log c}$, where $a=\log p/\log q$, and $c$ is as in Lemma
\ref{th4.5}. Since $\iota^*({\bf i}) =\sum_{j=1}^\infty\iota({\bf i}
\ 0^\infty|_j)q^{j-1}$ for any ${\bf i}\in\Sigma^\iota$.  Let $\ell$
be the largest index such that $\iota({\bf i} 0^\infty|_\ell)\ne 0$.
Then $\ell\ge n^b+1$,  and a direct estimation shows that
$$
|\iota^*({\bf i})|\ge
q^{\ell-1}-(q-2)\sum_{j=1}^{\ell-1}q^{j-1}\ge q^{\ell-3}+\frac12.
$$
 This together with the assumption on $\xi$ implies
 that
$$
\big |\rho^{m_0-1}\xi+\frac {\iota^*({\bf i})}{\rho N}\big | \ge
\frac{|\iota^*({\bf i} )|}{\rho N}-\frac1{2\rho N} \ge q^{\ell-4}\ge
q^{n^b-3}.
$$
Let  $\eta=\rho^{(n+1)^b}(\rho^{m_0-1}\xi+\frac {\iota^*({\bf
i})}{\rho N})$. It is easy to see that if $n$ large enough, then
$
(n+1)^b+3 \leq n^b+b^2n^{b-1}.
$
Hence if we take a large  $n_0$, then for $n\ge n_0$,
$$
|\eta|\ge \frac{p^{(n+1)^b}}{q^{(n+1)^b-n^b+3}}\ge
\frac{p^{(n+1)^b}}{q^{ b^2n^{b-1}}}
=\left(\frac{p^{(1+1/n)^bn}}{q^{b^2}}\right)^{n^{b-1}}\ge
\left(\frac{p^{n}}{q^{b^2}}\right)^{n^{b-1}}\ge q^{n^{b-1}}.
$$
Applying Lemma \ref{th4.5} to $\eta_i = \rho^i \eta$ recursively,
we have
$$
|\widehat{\mu}(\eta)|\le c|\widehat{\mu}(\eta_1)|\le \cdots\le
c^l|\widehat{\mu}(\eta_l)|
$$
as long as  $|\eta_l|\ge 1$. This is the case if we let $l=[ \log_a\frac
{2n^{1-b}}{1-a}]$, because
$$
|\eta_l|\ge \rho^2|\eta_{l-1}|^a\ge \cdots\ge
\rho^{2+2a+\cdots+2a^{l-1}}|\eta|^{a^l}\ge\rho^{\frac
2{1-a}}|\eta|^{a^l}>q^{-\frac
2{1-a}}q^{n^{b-1}a^l}\ge 1.
$$
Hence,
$$
|\widehat{\mu}(\eta)|\le c^l\le c^{\log_a \frac {2n^{1-b}}{1-a}}=c^{\log_a \frac 2{1-a}}\ {n^{-(b-1)\log c/\log a}}.
$$ The lemma
follows by assigning $C$ and $\beta$ in the obvious way. \eproof

\bigskip

{\noindent\bf Proof of Proposition \ref{prop3.2}.} We assume all the
parameters in Lemma \ref{th4.6}. To simplify the notations, we write
\ $\mu=
 \mu_{\rho, N}$,  $\alpha({\bf i}) =\rho^{-m_0}N^{-1}\iota^*({\bf i})$,\  ${\mathcal I}_{n} =\{{\bf i} \in \Omega_N^{\iota}, \, |{\bf
i}|\le n^b\}$,\  ${\mathcal I}_{n, n+1} =\{{\bf i} \in
\Omega_N^{\iota}, \, n^b< |{\bf i}|\le (n+1)^b\}$.  Let
$$
Q_n(\xi)={\sum}_{{\mathcal I}_n} |\widehat{\mu}\big(\xi+\alpha({\bf
i})\big)|^2.
$$
Then
\begin{eqnarray*}
Q_{n+1}(\xi) &=&Q_n(\xi)+{\sum}_{{\mathcal I}_{n, n+1}}
|\widehat{\mu}\big(\xi+\alpha ({\bf i})\big)|^2\\
&=&Q_n(\xi) +{\sum}_{{\mathcal I}_{n, n+1}}
|\widehat{\mu}_{(n+1)^b+m_0-1}\big(\xi+\alpha ({\bf i})\big)|^2
|\widehat{\mu}\left(\rho^{(n+1)^b+m_0-1}\big(\xi+
\alpha ({\bf i})\big)\right)|^2\\
&\le&Q_n(\xi)+\frac {C^2}{n^{2\beta}}{\sum}_{{\mathcal I}_{n, n+1}}
|\widehat{\mu}_{(n+1)^b+m_0-1}\big(\xi+\alpha ({\bf
i})\big)|^2  \qquad \hbox {(by Lemma \ref {th4.6})} \\
&\le &Q_n(\xi)+\frac {C^2}{n^{2\beta}}\left(1-{\sum}_{{\mathcal
I}_n} |\widehat{\mu}_{(n+1)^b+m_0-1}\big(\xi+\alpha ({\bf
i})\big)|^2\right) \quad\hbox {(by Lemma \ref {th4.5})}\\
&\le&Q_n(\xi)+\frac {C^2}{n^{2\beta}}\left(1-Q_n(\xi)\right).
\end{eqnarray*}
This implies that  $n > n_0$,
$$
1-Q_{n+1}(\xi)\ge\Big( 1-Q_{n}(\xi)\Big)\Big( 1-\frac
{C^2}{n^{2\beta}}\Big)\ge \cdots\ge \left(
1-Q_{n_0}(\xi)\right)\prod_{k=n_0}^n\Big ( 1-\frac {C^2}{k^{2\beta} }\Big ).
$$
Now let $ Q(\xi)={\sum}_{{\bf i} \in \Omega_N^{\iota}}
|\widehat{\mu}\big(\xi+\alpha({\bf i}\big)|^2$, it is the sum over a
maximal bi-zero set (by Proposition \ref{thm3.5}). The above implies
$$
1-Q(\xi)\ge C' \left( 1-Q_{n_0}(\xi)\right).
$$
where $C' =\prod_{k=n_0}^\infty( 1-\frac {C^2}{k^{2\beta}}) \not=0$.
This implies that $Q(\xi)\not \equiv 1$, and hence by Theorem
\ref{th2.1} and Proposition \ref{thm3.5}, any maximal bi-zero set of
$\mu_{\rho, N}$ cannot be a spectrum when $\rho =p/q, p, q$
co-prime, and $p\not=1$. \eproof

\bigskip
\bigskip
\section {\bf Remarks}

\bigskip

It was proved  in \cite {HLL} that if $\mu$ is a spectral
self-similar measure with support in $[0, 1]$ and $\nu$ is a
probability counting measure support on a finite set in $\Z$, then
the convolution $\mu\ast\nu$  is a spectral measure if and only if
$\nu$ is a spectral measure. It was pointed out by Gabardo and Lai
(private communication) that if both $\mu$ and $\nu$ are two
probability measures with $\mu\ast\nu=L|_{[0, 1]}$,  where
$L|_{[0,1]}$ is the Lebesgue measure restricted on $[0, 1]$, then
both $\mu$ and $\nu$ are spectral measures (which is a corollary of
the main results in \cite{AH} and \cite{L}).  It has been asked:

{\it Is the
convolution of two spectral self-similar measures with essentially
disjoint supports a spectral measure?}.

The  question can be answered by Theorem \ref{th1.1}. Observe that $\{0,1,2, 3\}
=\{0,1\} \oplus \{0,2\}$, hence
$$
\mu_{1/6, \ 4}=\mu_{1/6,\  2}\ast \mu_{1/6, \ \{0,  2\}}.
$$
It follows that both $\mu_{1/6,\ 2}$  and $\mu_{1/6,\ \{0,
2\}}$  are spectral measures (by \cite{JP} or Theorem \ref{th1.1}), but Theorem \ref{th1.1} implies that $\mu_{1/6, 4}$ is not a spectral measure. As a consequence, convolution of two spectral measures may not be spectral.

\bigskip

One of the challenge questions on the spectral measures is the conjecture of {\L}aba and Wang \cite{LW}:

{\it Let $\mu$ be a self-similar measure as in (\ref{eq1.1}), then
$\mu$  is a spectral measure if and only if (i) \ $w_j = 1/N$; (ii)
\ $\rho = 1/q$ for some integer $q>1$; and (iii)
there exist a constant $c$ and an integer digit set $\D'$ such that
$\D=c\D'$ and $ {\mathcal D}'\oplus {\mathcal B}\equiv \{0, \cdots,
q-1\} \ (\hbox {mod} \ q)$ \ for some ${\mathcal B} \subset {\Bbb
Z}$.}

In \cite{DL}, it was shown that (i) is necessary for a spectral
measure under  the no overlap condition. Our Theorem \ref{th1.1}
settles the case where ${\mathcal D }= \{0,\cdots, N-1\}$. The digit
set ${\mathcal D}'$ in (iii) is called an {\it integer tile}. The
study of integer tiles has a a long history related to the geometry
of numbers (\cite {CM} and the references there), and the spectral
property of ${\mathcal D}$ as a discrete set itself is still
unsolved \cite {L}.

\bigskip

As was proved in \cite {JP}, the Cantor measure $\mu_{1/k}$ with $k$ an odd integer is not a spectral measure. It is well known that a relaxing of the orthonormal basis is the concept of {\it frame} introduced by Duffin and Schaeffer in the 50's (see \cite {Chr}). We call a measure $\mu$  an {\it F-spectral measure} (F for frame) if there exists a countable set $\{e_\lambda: \lambda \in\Lambda\}$ and $A, B >0$ such that for any $f \in L^2(\mu)$,
$$
A ||f||^2 \leq \sum_{\lambda\in \Lambda} |\langle f, e_\lambda\rangle |^2 \leq B ||f||^2,
$$
and call $\mu$ a {\it R-spectral measure} if in addition it is a basis (R for Riesz). The frame structure of $L^2[0,1]$ has been studied in detail in \cite {Lan, OS};  also there are extensive studies of the frames on $L^2(\mu)$ \cite {DHL, P}. However the basic problem whether $\mu_{1/k}$ with $k$ an odd integer, in particular for $\mu_{1/3}$, is an F-(or R-)spectral measure is still unresolved.

\bigskip

{\it Acknowledgement:}  The authors like to thank Professor D.J. Feng and Dr. C.K. Lai for many helpful discussions.


\bigskip

\begin{thebibliography}{9999}
\smallskip

\bibitem {AH}
{\sc L.-X. An and X.-G. He}, {\it A class of spectral Moran
measures},  J. Funct. Anal. (To appear).

\bibitem  {CM}
{\sc E. Coven and A. Meyerowitz}, {\it Tiling the integers with translates of one finite set},
 J. Algebra   212(1999), 161--17.

 \bibitem {Chr}
{\sc O. Christensen}, {\it An Introduction to Frames and Riesz Bases}, Applied and Numerical Harmonic Analysis. Birkh\"{a}user Boston Inc., Boston, MA, 2003.

\bibitem {D}
{\sc X.-R. Dai}, {\it When does a Bernoulli convolution admit a spectrum?},  Adv. Math. 231(2012),  1681-1693.

\bibitem {DHL}
{\sc X.-R. Dai, X.-G. He and C.-K. Lai,} {\it Spectral property of
Cantor measures with consecutive digits},  Adv. Math. 242 (2013),
187-208.

\bibitem {Deng}
{\sc Q.-R. Deng}, {\it On a spectral property of self-similar
measures}, J. Math. Anal. Appl. (To appear).

\bibitem {DHS}
{\sc D. Dutkay, D. Han and Q. Sun}, {\it On spectra of a Cantor
measure}, Adv. Math. 221(2009), 251-276.


\bibitem {DHSW}
{\sc D. Dutkay, D. Han, Q. Sun and E. Weber}, {\it On the Beurling
dimension of exponential frames}, Adv. Math. 226(2011), 285-297.



\bibitem {DJ}
{\sc D. Dutkay and P. Jorgensen}, {\it Fourier frequencies in affine
iterated function systems}, J. Funct. Anal. 247(2007), 110-137.

\bibitem {DL}
{\sc D. Dutkay and C.-K. Lai}, {\it Uniformity of measures with
Fourier frames}, preprint.


\bibitem {Fal}
{\sc K. J. Falconer}, {\it Fractal Geometry, Mathematical
Foundations and Applications}, Wiley, New York, 1990.

\bibitem {F}
{\sc B. Fuglede}, {\it Commuting self-adjoint partial differential
operators and a group theoretic problem}, J. Funct. Anal. 16(1974), 101-121.

\bibitem {HLL}
{\sc X.-G. He, C.-K. Lai and K.-S. Lau}, {\it Exponential spectra in $L^2(\mu)$}, Appl. Comput. Harmon. Anal.  34(2013) 327-338.


\bibitem  {HuL}
{\sc T.-Y. Hu and K.-S. Lau}, {\it Spectral property of the
Bernoulli convolutions}, Adv. Math. 219(2008), 554-567.



\bibitem  {JP}
{\sc P. Jorgensen and S. Pedersen}, {\it Dense analytic subspaces in
fractal $L^2$ spaces}, J. Anal. Math. 75(1998), 185-228.

\bibitem {KM}
{\sc M. Kolountzakis and M. Matolcsi}, {\it  Tiles with no spectra},
Forum Math. 18(2006), 519-528.


\bibitem {L}
{\sc I. $\L$aba,} {\it The spectral set conjecture and multiplicative properties of roots of polynomials},
\newblock J. London Math. Soc.  65(2002), 661--671.

\bibitem  {LW}
{\sc I. {\L}aba and Y. Wang}, {\it On spectral Cantor measures}, J.
Funct. Anal.  193(2002), 409-420.


\bibitem {Lai}
{\sc C.-K. Lai}, {\it On Fourier frame of absolutely continuous
measures}, J. Funct. Anal. 261(2011), 2877-2889.

\bibitem  {Lan}
{\sc H. Landau}, {\it Necessary density conditions for sampling and
interpolation of certain entire functions}, Acta Math. 117(1967), 37-52.

\bibitem{L}
{\sc T. Lewis}, {\it The factorization of the rectangular
distribution}, J. Applied Probab.  4(1967), 529-542.


\bibitem {Li1}
{\sc J.-L. Li}, {\it $\mu_{M,D}-$orthogonality and compatible pair},
J. Funct. Anal.  244(2007), 628-638.

\bibitem {Li2}
{\sc J.-L. Li}, {\it Spectra of a class of self-affine measures}, J.
Funct. Anal. 260(2011),  1086-1095.


\bibitem  {OS}
{\sc J. Ortega-Cerd\`{a} and K. Seip}, {\it Fourier frames}, Annal Math.  155(2002), 789-806.

\bibitem {P}
{\sc A. Poltoratski}, {\it A problem on completeness of
exponentials}, Annal Math. (To appear).


 \bibitem  {S}
 {\sc R. Strichartz}, {\it Convergence of Mock Fourier series }, J. Anal.
 Math. 99(2006), 333-353.

 \bibitem  {T}
  {\sc T. Tao}, {\it Fuglede's conjecture is false in 5 or higher dimensions}, Math. Res. Lett.  11(2004), 251-258.


\end{thebibliography}
\end{document}